\documentclass[12pt,a4paper,reqno]{amsart}
\usepackage{amssymb,enumitem,mathtools}\usepackage[a4paper,margin=1in]{geometry}
\usepackage{amsmath,amssymb,amsthm,mathtools}
\usepackage{mathrsfs}
\usepackage{bm}
\usepackage{xcolor}
\usepackage{enumitem}
\usepackage{hyperref}
\usepackage{microtype}
\usepackage{booktabs}
\usepackage{lineno}
\usepackage{graphicx}
\allowdisplaybreaks
\usepackage{hyperref}
\hypersetup{
colorlinks=true,
urlcolor=black,
citecolor=blue}
\numberwithin{equation}{section}
\usepackage{color}

\hypersetup{
  colorlinks=true,
  linkcolor=blue!55!black,
  citecolor=blue!55!black,
  urlcolor=blue!55!black
}
\usepackage{amsthm}

\theoremstyle{plain}

\newtheorem{theorem}{Theorem}[section]
\newtheorem{proposition}[theorem]{Proposition}
\newtheorem{lemma}[theorem]{Lemma}
\newtheorem{corollary}[theorem]{Corollary}

\newtheorem{example}[theorem]{Example}
\theoremstyle{definition}

\theoremstyle{remark}
\newtheorem{remark}[theorem]{Remark}

\usepackage{amsthm}

\theoremstyle{plain}
\newtheorem{theoremA}{Theorem}

\newtheorem{theoremB}{Theorem}

\newtheorem{corollaryC}{Corollary}

\newtheorem{question}[theorem]{Question}

\newcommand{\tr}{\operatorname{tr}}
\newcommand{\A}{\mathcal A}
\newcommand{\Rcheck}{\check{R}}
\newcommand{\Rm}{\mathrm{Rm}}
\newcommand{\Ric}{\operatorname{Ric}}
\newcommand{\Scal}{\operatorname{s}}
\newcommand{\grad}{\operatorname{grad}}
\newcommand{\Hess}{\operatorname{Hess}}

\title[Rigidity for the quadratic Riemann curvature functional]
{Rigidity for the quadratic Riemann curvature functional in dimensions greater than four}
\author{Yunhee Euh}
\address{Department of Mathematics, Sungkyunkwan University, Suwon 16419, Korea}
\email{prettyfish@skku.edu}

\author{JeongHyeong Park}
\address{Department of Mathematics, Sungkyunkwan University, Suwon 16419, Korea}
\email{parkj@skku.edu}


\begin{document}

\begin{abstract}
We investigate rigidity properties of critical metrics for the quadratic Riemann curvature functional on Riemannian manifolds of dimension greater than four. This functional is the integral of the squared norm of the Riemann curvature tensor. On noncompact manifolds, criticality is understood with respect to compactly supported variations, without assuming finite total curvature energy. We prove that every complete noncompact critical metric with finite curvature energy is either flat or has strictly negative scalar curvature everywhere. We further show that every complete critical metric with concircular scalar curvature is flat. Here concircularity means that the Hessian of the scalar curvature is pointwise proportional to the metric. Consequently, every complete weakly Einstein critical metric with divergence-free Cotton tensor is flat. The latter two flatness results require neither noncompactness nor finite curvature energy. Finally, we construct a complete nonflat critical metric in dimension five with vanishing Cotton tensor, showing that the concircularity and weakly Einstein assumptions are essential in the respective flatness results.

\end{abstract}

\maketitle

\section{Introduction}\label{sec:1_intro}
Quadratic curvature functionals are natural variational objects in
Riemannian geometry. Among them, the quadratic Riemann curvature
functional
\begin{equation}\label{eq:Afunctional}
  \mathcal A(g)=\int_M |\Rm|^2\,dv_g
\end{equation}
is distinguished by the fact that $\mathcal A(g)=0$ precisely for flat
metrics. We say that $g$ has finite curvature energy if $\mathcal A(g)<\infty$.
Berger initiated the systematic variational study of quadratic
curvature functionals \cite{Berger1970}, and the basic first-variation
formulas are recorded, for example, in Besse's monograph
\cite{Besse1987}.

Critical metrics of quadratic curvature functionals have been studied
on both closed and complete noncompact manifolds. Catino obtained
pointwise and integral Weitzenb\"ock formulas for critical metrics of
general quadratic curvature functionals on closed manifolds
\cite{Catino2015}. More recently, Catino--Mastrolia--Monticelli
established flatness and scalar-flatness results for complete
noncompact critical metrics of quadratic functionals involving
$|\Ric|^2$ and $\Scal^2$, where $\Ric$ and $\Scal$ denote the Ricci tensor
and the scalar curvature, respectively \cite{CMM2023}. For the quadratic scalar curvature functional, scalar-flatness of
complete noncompact critical metrics satisfying $\Scal\in L^2(M)$ was established in dimensions $n\ge10$ in \cite{CMM2026} and in dimensions $5\le n\le9$ in \cite{Zhang2026}. These results motivate the study of the {{quadratic}} Riemann curvature functional \eqref{eq:Afunctional} on complete noncompact manifolds. In this setting, the traced Euler–Lagrange equation implies that the scalar curvature is subharmonic when $n>4$. We exploit this property to obtain restrictions on the scalar curvature and flatness under additional geometric assumptions.

We investigate this problem, as well as related rigidity questions,
in connection with the weakly Einstein condition
\begin{equation*}\label{eq:WEintro}
  \Rcheck=\frac{|\Rm|^2}{n}g,
\end{equation*}
where $\Rcheck_{ij}=R_{iabc}R_j{}^{abc}$. Here and throughout the paper, repeated indices are summed
according to the Einstein summation convention.
The weakly Einstein condition was introduced in dimension four, as a generalization of Einstein spaces, defining metrics for which a curvature identity valid for Einstein four-manifolds continues to hold \cite{EPS2013a}, and was studied further in \cite{EPS2013b}.
It also arises naturally in the Euler--Lagrange equation for
$\mathcal A$. In dimensions greater than four, neither the weakly
Einstein condition nor the Einstein condition implies the other in
general. We establish rigidity results for critical metrics of $\mathcal A$, both with and without the weakly Einstein assumption.

The properties of critical metrics for $\mathcal A$
depend significantly on the dimension.
In dimensions two and three, known results for
quadratic curvature functionals imply that
complete critical metrics for $\mathcal A$
are flat \cite{CMM2023}.
Dimension four is exceptional: the trace of the
Euler--Lagrange equation for $\mathcal A$
reduces to $\Delta\Scal=0$, and nonflat examples exist.
In particular, the product metric on the closed manifold
$S^2(1)\times\Sigma(-1)$ is weakly Einstein and
critical for $\mathcal A$, but not Einstein.
Here $S^2(1)$ denotes the unit round $2$-sphere,
and $\Sigma(-1)$ is a closed surface of constant
Gaussian curvature $-1$
\cite{EuhParkSekigawa2011}. {{In our previous work \cite{EuhPark2025}, we studied critical metrics for
$\mathcal A$ with finite curvature energy on complete
four-dimensional manifolds.}}

For $n>4$, the trace of the
Euler--Lagrange equation for $\mathcal A$ on the space of all metrics
yields
\begin{equation}\label{eq:trace-intro}
  \Delta\Scal=\frac{n-4}{4}|\Rm|^2\geq0.
\end{equation}
Thus the scalar curvature is subharmonic, leading to a different
rigidity problem in higher dimensions.

We ask the following question.

\begin{question}\label{ques:finite-energy-rigidity}
Let $(M^n,g)$, $n>4$, be a complete noncompact weakly Einstein
Riemannian manifold. Suppose that $g$ is a critical metric for
$\mathcal A$ on the space $\mathcal M(M)$ of all Riemannian metrics
on $M$ and has finite curvature energy. Does it follow that $(M,g)$ is flat?
\end{question}

We do not completely resolve
Question~\ref{ques:finite-energy-rigidity}.
Our first main result restricts the scalar curvature
of complete noncompact metrics that are critical
for $\mathcal A$ and have finite curvature energy, without
requiring the weakly Einstein condition.
\begin{theoremA}\label{thm:theoremA}
Let $(M^n,g)$, $n>4$, be a complete noncompact Riemannian manifold.
Suppose that $g$ is a critical metric for $\mathcal A$ on
$\mathcal M(M)$ and has finite curvature energy. Then either $(M,g)$
is flat or its scalar curvature $\Scal$ is strictly negative everywhere
on $M$.
\end{theoremA}

Thus a nonflat metric satisfying the hypotheses of
Question~\ref{ques:finite-energy-rigidity} must have strictly negative
scalar curvature everywhere. Theorem~\ref{thm:theoremA} follows from
\eqref{eq:trace-intro}, a finite-energy cutoff argument, and the
maximum principle.

Our second main result establishes flatness when the scalar curvature is concircular, without requiring either the weakly Einstein condition or finite curvature energy. Recall that a
smooth function $f$ is concircular if {{its Hessian is pointwise proportional to the metric, that is}}
$$\Hess f=\varphi g$$ for some
smooth function $\varphi$, or equivalently, if
$\Hess f=(\Delta f/n)g$.

\begin{theoremB}\label{thm:concircular}
Let $(M^n,g)$, $n>4$, be a complete Riemannian manifold. Suppose
that $g$ is a critical metric for $\mathcal A$ on $\mathcal M(M)$.
If its scalar curvature is concircular, then $(M,g)$ is flat.
\end{theoremB}

For a weakly Einstein metric that is critical
for $\mathcal A$, we establish the differential identity
\begin{equation}\label{eq:divC-intro}
  \nabla^a C_{ija}
  =\frac{1}{2(n-1)}
  \left((\Hess\Scal)_{ij}-\frac{\Delta\Scal}{n}g_{ij}\right),
\end{equation}
where $C$ denotes the Cotton tensor; see {Proposition~\ref{prop:divergence-identities}}.
Hence, the Cotton tensor is
divergence-free if and only if the scalar curvature is concircular.
Combining \eqref{eq:divC-intro} with Theorem~\ref{thm:concircular} gives the following consequence.

\begin{corollaryC}\label{cor:divC-rigidity-intro}
Let $(M^n,g)$, $n>4$, be a complete weakly Einstein Riemannian
manifold. Suppose that $g$ is a critical metric for $\mathcal A$ on
$\mathcal M(M)$. If its Cotton tensor is divergence-free, then
$(M,g)$ is flat.
\end{corollaryC}

In particular, the conclusion holds if the Cotton tensor vanishes.
Corollary~\ref{cor:divC-rigidity-intro} requires neither noncompactness
nor finite curvature energy and therefore gives an affirmative answer
to Question~\ref{ques:finite-energy-rigidity} under the additional
assumption that the Cotton tensor is divergence-free.

We also study weakly Einstein metrics on closed manifolds that are
critical for $\mathcal A$ among unit-volume metrics. In this setting,
we derive an integral identity relating the Cotton tensor to the
gradient of the scalar curvature.

The paper is organized as follows. In Section~\ref{sec:2_notation}, we fix the notation
and recall the weakly Einstein condition and the tensors used
throughout the paper. In Section~\ref{sec:3_EL}, we review the Euler--Lagrange
equations for $\mathcal A$ on $\mathcal M(M)$ and on the space
$\mathcal M_1(M)$ of unit-volume Riemannian metrics. In Section~\ref{sec:4_E},
we derive integral identities for closed weakly Einstein metrics that
are critical for $\mathcal A$ on $\mathcal M_1(M)$. In Section~\ref{sec:5_TheoremA}, we study complete noncompact metrics
that are critical for $\mathcal A$ and have finite curvature 
energy, and prove Theorem~\ref{thm:theoremA}.
In Section~\ref{sec:6_identities}, we establish differential identities
for weakly Einstein metrics that are critical
for $\mathcal A$, including \eqref{eq:divC-intro}. Finally, in Section~\ref{sec:7_TheoremB}, we prove Theorem~\ref{thm:concircular} and
Corollary~\ref{cor:divC-rigidity-intro}. We also construct a complete nonflat metric in dimension five that is critical for $\mathcal A$ and has vanishing Cotton tensor. This example shows that the concircularity assumption in Theorem B and the weakly Einstein assumption in Corollary C cannot be omitted.

\section{Preliminaries and Conventions}\label{sec:2_notation}

All manifolds are assumed to be connected and without boundary unless otherwise stated. We use the curvature convention
\begin{equation}\label{eq:curvature-convention}
  R(X,Y)Z=[\nabla_X,\nabla_Y]Z-\nabla_{[X,Y]}Z
\end{equation}
and we set
$  \Rm(X,Y,Z,W)=g(R(X,Y)Z,W)$ for $X,Y,Z,W\in\mathfrak{X}(M)$, where $\mathfrak{X}(M)$ denotes the space of smooth vector fields on $M$.
Thus, with respect to a local orthonormal frame $\{e_i\}$,
$  R_{ijkl}=\Rm(e_i,e_j,e_k,e_l).$
The Ricci tensor $\Ric$ and scalar curvature $\Scal$ are defined by
\begin{equation*}\label{eq:ricci-scalar-convention}
  \Ric(X,Y)=\sum_{a=1}^n \Rm(e_a,X,Y,e_a),
  \qquad
  \Scal=\tr_g\Ric.
\end{equation*}

{{For any symmetric $(0,2)$-tensor $T$, we denote its trace-free part by
\begin{equation*}
T^\circ_{ij}
:=
T_{ij}-\frac{\tr_g T}{n}g_{ij}.
\end{equation*}}}
In particular, the trace-free Ricci tensor is given by
\begin{equation*}\label{eq:Edef}
  \Ric^{\circ}_{ij}=\Ric_{ij}-\frac{\Scal}{n}g_{ij}.
\end{equation*}
The contracted Bianchi identity gives
\begin{equation}\label{eq:divE}
  \nabla^i\Ric^{\circ}_{ij}=\frac{n-2}{2n}\nabla_j\Scal.
\end{equation}
We also have
\begin{equation}\label{eq:trac-free}
\begin{gathered}
\Delta \Ric^{\circ}_{ij}
  =
  (\Delta\Ric)^\circ_{ij}
  =
  \Delta\Ric_{ij}
  -\frac{\Delta\Scal}{n}g_{ij},\\
  (\Hess\Scal)^\circ_{ij}
 =
  (\Hess\Scal)_{ij}
  -\frac{\Delta\Scal}{n}g_{ij},
\end{gathered}
\end{equation}
where $\Delta=g^{ab}\nabla_a\nabla_b$ denotes the Laplace--Beltrami operator on functions and the rough Laplacian on tensor fields. 

We define the symmetric 
$(0,2)$-tensor $\check R$ by
\begin{equation*}\label{eq:Rcheckdef}
  \Rcheck_{ij}=R_{iabc}R_{j}{}^{abc}.
\end{equation*}
Since $\tr_g\Rcheck=|\Rm|^2$, the weakly Einstein condition is
\begin{equation}\label{eq:WE}
  \Rcheck_{ij}=\frac{|\Rm|^2}{n}g_{ij},
\end{equation}
that is, $\Rcheck^\circ=0$.

The Cotton tensor is defined by
\begin{equation*}\label{eq:Cotton}
  C_{ijk}
  =\nabla_k\Ric_{ij}-\nabla_j\Ric_{ik}
  -\frac{1}{2(n-1)}
   \bigl(\nabla_k\Scal\,g_{ij}-\nabla_j\Scal\,g_{ik}\bigr).
\end{equation*}
The divergence of the Cotton tensor
is taken in its third index:
\begin{equation*}
  (\operatorname{div}C)_{ij}:=\nabla^a C_{ija}.
\end{equation*}
We say that $C$ is divergence-free if $\nabla^a C_{ija}=0$.

\section{Euler--Lagrange Equations and Criticality}\label{sec:3_EL}
With the convention \eqref{eq:curvature-convention}, the formal
$L^2$-gradient of $\A$ is given by Berger's formula
\cite{Berger1970,Besse1987,EuhPark2025}:
\begin{equation}\label{eq:gradA}
\begin{split}
  (\grad\A)_{ij}
  ={}&-2R_{iabc}R_j{}^{abc}
      -4\Delta\Ric_{ij}
      +2(\Hess\Scal)_{ij}\\
      &      +4\Ric_{ia}\Ric_j{}^a
      -4\Ric^{ab}R_{iabj}
      +\frac12|\Rm|^2g_{ij}.
\end{split}
\end{equation}

On a noncompact manifold $M$, we define criticality
using compactly supported variations, without assuming that
$\A(g)$ is finite. Given a smooth compactly supported symmetric
$(0,2)$-tensor $h$, set $g_t=g+th$ for sufficiently small $|t|$,
and choose a relatively compact domain $\Omega$ with smooth
boundary such that $\operatorname{supp}h\subset\Omega$.
We say that $g$ is critical for $\A$ on $\mathcal M(M)$ if
\begin{equation*}
  \left.\frac{d}{dt}\right|_{t=0}
  \int_\Omega |\Rm_{g_t}|_{g_t}^2\,dv_{g_t}=0
\end{equation*}
for every such $h$. Since $g_t=g$ outside
$\operatorname{supp}h$, this definition is independent of
the choice of $\Omega$. The local first-variation formula gives
\begin{equation*}
  \left.\frac{d}{dt}\right|_{t=0}
  \int_\Omega |\Rm_{g_t}|_{g_t}^2\,dv_{g_t}
  =
  \int_\Omega \langle\grad\A,h\rangle_g\,dv_g.
\end{equation*}
Thus criticality is equivalent to
\begin{equation*}
  \grad\A=0.
\end{equation*}
For a closed manifold, if $g$ is critical for $\A$ on $\mathcal{M}_1(M)$, then
\begin{equation*}\label{eq:critical-M1}
  \grad\A=cg
\end{equation*}
for a constant Lagrange multiplier $c$. Hence the trace-free part of the Euler--Lagrange equation is the same in both cases.
\medskip

Suppose first that $g$ is critical for $\mathcal A$ on $\mathcal M(M)$.
Then $\grad\A=0$. Taking the trace of \eqref{eq:gradA} 
we obtain
\begin{equation}\label{eq:trace-equation}
  \Delta\Scal=\frac{n-4}{4}|\Rm|^2.
\end{equation}
This identity is the basic reason that dimension four is exceptional.

For comparison, suppose that $M$ is closed and that $g$ is critical for
$\mathcal A$ on $\mathcal M_1(M)$. Then
$
  \grad\A=cg
$
for some constant $c$. Taking the trace of the Euler--Lagrange equation gives
$$  \Delta\Scal  =  \frac{n-4}{4}|\Rm|^2-\frac n2 c.$$
Since $g\in\mathcal{M}_1(M)$, integration over $M$ yields
$$  c=\frac{n-4}{2n}\mathcal A(g).$$
Consequently,
\begin{equation}\label{eq:trace-unit-volume}
  \Delta\Scal
  =
  \frac{n-4}{4}
  \left(|\Rm|^2-\mathcal A(g)\right).
\end{equation}
Although the trace equations differ in the two variational problems,
both Euler--Lagrange equations imply
$(\grad\A)^\circ=0$. In the next section, we use this common trace-free equation
to derive a differential equation for the trace-free Ricci tensor
of weakly Einstein metrics that are critical for $\mathcal A$.

\section{The Trace-Free Ricci Tensor and Closed Critical Metrics}\label{sec:4_E}
We first derive a differential equation for the trace-free Ricci tensor of weakly Einstein critical metrics.

\begin{lemma}\label{thm:local-main}
Let $(M^n,g)$, $n>4$, be a weakly Einstein Riemannian manifold. Suppose that either $g$ is a critical metric for $\mathcal A$ on  $\mathcal M(M)$, or $M$ is closed and $g$ is a critical metric for $\mathcal A$ on  $\mathcal M_1(M)$. Then $\Ric^{\circ}$ satisfies
\begin{equation}\label{eq:E-equation-intro}
\begin{split}
  \Delta \Ric^{\circ}_{ij}
  ={}&
  \frac12(\Hess\Scal)_{ij}
  -\frac{1}{2n} \Delta\Scal\, g_{ij}
  \\ &{}-(\Ric^{\circ})^{ab}R_{iabj}
  +\Ric^{\circ}_{ia}\Ric^{\circ}_j{}^a
  +\frac{\Scal}{n}\Ric^{\circ}_{ij}.
\end{split}
\end{equation}
Consequently,
\begin{equation}\label{eq:bochner-intro}
\begin{split}
  \frac12\Delta |\Ric^{\circ}|^2
  ={}&
  |\nabla \Ric^{\circ}|^2
  +\frac12 (\Ric^{\circ})^{ab}(\Hess\Scal)_{ab}
  \\ &{}-(\Ric^{\circ})^{ab}R_{cabd}(\Ric^{\circ})^{cd} 
  +(\Ric^{\circ})^{ab}\Ric^{\circ}_{ac}\Ric^{\circ}_b{}^c
  +\frac{\Scal}{n}|\Ric^{\circ}|^2.
\end{split}
\end{equation}
\end{lemma}

\begin{proof}
In either variational problem, the discussion in Section~\ref{sec:3_EL} implies that
$(\grad\A)^\circ=0$.
Since $g$ is weakly Einstein, \eqref{eq:WE} gives
$\Rcheck^\circ=0$. Hence, taking the trace-free part of the
Euler--Lagrange equation yields
\begin{equation}\label{eq:tracefree-grad}
  0
  =
  -4(\Delta\Ric)^\circ_{ij}
  +2(\Hess\Scal)^\circ_{ij}
  +4\bigl(\Ric_{ia}\Ric_j{}^a\bigr)^\circ
  -4\bigl(\Ric^{ab}R_{iabj}\bigr)^\circ.
\end{equation}
Using $\Ric=\Ric^{\circ}+\frac{\Scal}{n}g$, we compute
\begin{equation*}
\begin{aligned}
  \Ric_{ia}\Ric_j{}^a
  &=
\Ric^{\circ}_{ia}\Ric^{\circ}_j{}^a
  +\frac{2\Scal}{n}\Ric^{\circ}_{ij}
  +\frac{\Scal^2}{n^2}g_{ij},\\
  \Ric^{ab}R_{iabj}
  &=
  (\Ric^{\circ})^{ab}R_{iabj}
  +\frac{\Scal}{n}\Ric^{\circ}_{ij}
  +\frac{\Scal^2}{n^2}g_{ij}.
\end{aligned}
\end{equation*}
Therefore,
\begin{equation*}
  \Ric_{ia}\Ric_j{}^a-\Ric^{ab}R_{iabj}
  =
  \Ric^{\circ}_{ia}\Ric^{\circ}_j{}^a
  -(\Ric^{\circ})^{ab}R_{iabj}
  +\frac{\Scal}{n}\Ric^{\circ}_{ij}.
\end{equation*}
Since both sides are trace-free,
taking the trace-free part of this identity gives
\begin{equation}\label{eq:rhorho}
  \bigl(\Ric_{ia}\Ric_j{}^a\bigr)^\circ
  -
  \bigl(\Ric^{ab}R_{iabj}\bigr)^\circ
  =
  \Ric^{\circ}_{ia}\Ric^{\circ}_j{}^a
  -(\Ric^{\circ})^{ab}R_{iabj}
  +\frac{\Scal}{n}\Ric^{\circ}_{ij}.
\end{equation}
Using \eqref{eq:trac-free} to express
$(\Delta\Ric)^\circ$ and $(\Hess\Scal)^\circ$ in
\eqref{eq:tracefree-grad}, together with \eqref{eq:rhorho}, we obtain
$$
  \Delta \Ric^{\circ}_{ij}
  =
  \frac12
  \left(
    (\Hess\Scal)_{ij}
    -\frac{\Delta\Scal}{n}g_{ij}
  \right)
  -(\Ric^{\circ})^{ab}R_{iabj}
  +\Ric^{\circ}_{ia}\Ric^{\circ}_j{}^a
  +\frac{\Scal}{n}\Ric^{\circ}_{ij},
$$
which is \eqref{eq:E-equation-intro}.
Taking the inner product of \eqref{eq:E-equation-intro} with $\Ric^{\circ}$ and using
$$
  \frac12\Delta|\Ric^{\circ}|^2
  =
  |\nabla \Ric^{\circ}|^2+\langle\Delta \Ric^{\circ},\Ric^{\circ}\rangle,
$$
we obtain \eqref{eq:bochner-intro}.
\end{proof}

\begin{theorem}\label{thm:closed-main}
Let $(M^n,g)$, $n>4$, be a closed weakly Einstein Riemannian manifold. Suppose that $g$ is a critical metric for $\mathcal A$ on $\mathcal M_1(M)$. Then
\begin{equation*}\label{eq:Cotton-identity-intro}
  \int_M |C|^2\,dv_g
  =
  \frac{n-2}{2n(n-1)}
  \int_M |\nabla\Scal|^2\,dv_g.
\end{equation*}
 In particular, the Cotton tensor vanishes identically if and only if the scalar curvature is constant.
\end{theorem}

\begin{proof}
Catino proved an integral identity for the trace-free Ricci tensor \cite[Proposition~4.1]{Catino2015}.  
Rewritten in our curvature convention  \eqref{eq:curvature-convention},
Catino's identity takes the following form:
\begin{equation}\label{eq:CatinoP41}
\begin{aligned}
  &\int_M\left(
    |\nabla \Ric^{\circ}|^2
    -\frac{(n-2)^2}{4n(n-1)}|\nabla\Scal|^2
    -\frac12|C|^2
  \right)\,dv_g\\
  &={}\int_M\Bigl(
    (\Ric^{\circ})^{ab}R_{iabj}(\Ric^{\circ})^{ij}-(\Ric^{\circ})^{ij}\Ric^{\circ}_{ik}\Ric^{\circ}_{j}{}^k
    -\frac{\Scal}{n}|\Ric^{\circ}|^2
  \Bigr)\,dv_g.
\end{aligned}
\end{equation}

On a closed manifold, integrating \eqref{eq:bochner-intro} and using \eqref{eq:divE} gives
\begin{equation}\label{eq:integrated-bochner}
\begin{aligned}
  &\int_M\left(
    |\nabla \Ric^{\circ}|^2-\frac{n-2}{4n}|\nabla\Scal|^2
  \right)\,dv_g\\
  &={}\int_M\Bigl(
    (\Ric^{\circ})^{ab}R_{iabj}(\Ric^{\circ})^{ij}-(\Ric^{\circ})^{ij}\Ric^{\circ}_{ik}\Ric^{\circ}_{j}{}^k
    -\frac{\Scal}{n}|\Ric^{\circ}|^2
  \Bigr)\,dv_g,
\end{aligned}
\end{equation}
where we used
$$
  \frac12\int_M (\Ric^{\circ})^{ij}(\Hess\Scal)_{ij}\,dv_g
  =-\frac{n-2}{4n}\int_M|\nabla\Scal|^2\,dv_g.
$$
Comparing \eqref{eq:CatinoP41} and \eqref{eq:integrated-bochner}, we observe that their right-hand sides are identical. Subtracting the two identities yields
$$  \frac12\int_M|C|^2\,dv_g
  =\left[
    \frac{n-2}{4n}
    -\frac{(n-2)^2}{4n(n-1)}
   \right]
   \int_M|\nabla\Scal|^2\,dv_g.$$
We obtain
$$  \int_M|C|^2\,dv_g
  =\frac{n-2}{2n(n-1)}\int_M|\nabla\Scal|^2\,dv_g.$$
Since the coefficient $\frac{n-2}{2n(n-1)}$ is positive, $C\equiv0$ if and only if $\nabla\Scal\equiv0$. Since $M$ is connected, this is equivalent to $\Scal$ being constant.
\end{proof}

\begin{corollary}\label{cor}
Let $(M^n,g)$, $n>4$, be a closed weakly Einstein Riemannian manifold. Suppose that $g$ is a critical metric for $\mathcal A$ on $\mathcal M_1(M)$. Then the following conditions are equivalent:
\begin{enumerate}
\item The Cotton tensor $C$ vanishes identically;
\item The scalar curvature $\Scal$ is constant;
\item The squared norm $|\Rm|^2$ of the Riemann curvature tensor is constant.
\end{enumerate}
\end{corollary}

\begin{proof}
By Theorem~\ref{thm:closed-main}, the Cotton tensor $C$ vanishes
identically if and only if the scalar curvature $\Scal$ is constant.

Suppose that $\Scal$ is constant. Then $\Delta\Scal=0$, and
\eqref{eq:trace-unit-volume}, together with $n>4$, gives
$$
  |\Rm|^2=\mathcal A(g).
$$
Hence $|\Rm|^2$ is constant.

Conversely, suppose that $|\Rm|^2$ is constant. Since
$\operatorname{Vol}(M,g)=1$, we have
$$
  \mathcal A(g)
  =
  \int_M |\Rm|^2\,dv_g
  =
  |\Rm|^2.
$$
It follows from \eqref{eq:trace-unit-volume} that $\Delta\Scal=0$.
Since $M$ is closed, $\Scal$ must be constant.

Therefore, the three conditions are equivalent.
\end{proof}

When $g$ is a critical metric for $\mathcal A$ on  $\mathcal M(M)$, the assumption that $M$ is closed yields the following stronger rigidity result.

\begin{proposition}
Let $(M^n,g)$, $n>4$, be a closed Riemannian manifold. Suppose that $g$ is a critical metric for $\mathcal A$ on $\mathcal M(M)$. Then $(M,g)$ is flat.
\end{proposition}
\begin{proof}
{{Since $g$ is critical for $\mathcal A$ on $\mathcal M(M)$,
we have
\begin{equation*}
  \Delta\Scal=\frac{n-4}{4}|\Rm|^2.
\end{equation*}
Integrating over the closed manifold $M$ yields
\begin{equation*}
  0=\int_M\Delta\Scal\,dv_g
  =\frac{n-4}{4}\int_M|\Rm|^2\,dv_g.
\end{equation*}}}
Since $n>4$, it follows that $\Rm\equiv0$.
\end{proof}
Thus, in the closed case, criticality for $\mathcal A$ on
$\mathcal M(M)$ forces flatness when $n>4$.
For complete noncompact manifolds, however, the preceding integration
argument is no longer available. In the next section, we combine the
subharmonicity of the scalar curvature with the finite-energy assumption
to study the corresponding rigidity problem.

\section{Proof of Theorem~\ref{thm:theoremA} }\label{sec:5_TheoremA}
We now turn to complete noncompact metrics that are
critical for $\mathcal A$ and have finite curvature energy. Unlike the closed case, the traced Euler--Lagrange equation cannot be integrated directly to obtain flatness. Instead, we use the finite-energy assumption to control the scalar curvature and apply a cutoff argument to its positive part. We begin with the following standard curvature inequalities.

\begin{lemma}\textup{(\cite[p.~74]{BGM})}\label{lem:tau-L2}
On any Riemannian manifold of dimension $n$,
\begin{equation*}\label{eq:ineq}
  \Scal^2\leq n|\Ric|^2
  \qquad\text{and}\qquad
  |\Ric|^2\leq \frac{n-1}{2}|\Rm|^2.
\end{equation*}
Consequently, if $g$ has finite curvature energy,
or equivalently, $\Rm\in L^2(M)$, then both
$\Ric$ and $\Scal$ belong to $L^2(M)$.
\end{lemma}

We now prove Theorem~\ref{thm:theoremA}.

\begin{proof}[Proof of Theorem~\ref{thm:theoremA}]
By \eqref{eq:trace-equation} and $n>4$, we have
$\Delta\Scal\geq0$, so $\Scal$ is subharmonic.
Let
\begin{equation*}
  \Scal_+:=\max\{\Scal,0\}
\end{equation*}
denote the positive part of $\Scal$.
Since $\Scal$ is subharmonic, $\Scal_+$ is also
subharmonic in the distributional sense.
Moreover, Lemma~\ref{lem:tau-L2} implies that
$\Scal_+\in L^2(M)$.

Fix a point $p\in M$, and let $B_r=B_r(p)$
denote the geodesic ball of radius $r$ centered at $p$.
Choose smooth compactly supported cutoff functions
$\eta_r\in C_c^\infty(M)$ satisfying
\begin{equation*}
\begin{gathered}
  0\leq\eta_r\leq1,\qquad
  \eta_r\equiv1\ \text{on }B_r,\qquad
  \operatorname{supp}\eta_r\subset B_{2r},\qquad
  |\nabla\eta_r|\leq\frac{c_0}{r},
\end{gathered}
\end{equation*}
where $c_0>0$ is a constant independent of $r$.
Using the standard approximation argument,
we may test the weak subharmonicity inequality
for $\Scal_+$ with $\eta_r^2\Scal_+$.
This is the usual cutoff argument underlying
Yau's $L^p$-Liouville theorem
\cite{Yau1976}. 
By Young's inequality,
\begin{equation*}
\begin{aligned}
\int_M \eta_r^2|\nabla\Scal_+|^2\,dv_g
&\leq
2\int_M \eta_r\Scal_+|\nabla\eta_r|
|\nabla\Scal_+|\,dv_g\\
&\leq
\frac12\int_M \eta_r^2|\nabla\Scal_+|^2\,dv_g
+2\int_M \Scal_+^2|\nabla\eta_r|^2\,dv_g.
\end{aligned}
\end{equation*}
Consequently,
\begin{equation*}
\int_M \eta_r^2|\nabla\Scal_+|^2\,dv_g
\leq
4\int_M \Scal_+^2|\nabla\eta_r|^2\,dv_g
\leq
\frac{4c_0^2}{r^2}\int_M \Scal_+^2\,dv_g.
\end{equation*}
Since $\Scal_+\in L^2(M)$, the right-hand side
tends to zero as $r\to\infty$.
Thus $\nabla\Scal_+=0$ almost everywhere on $M$.
Since $\Scal_+$ is locally Lipschitz and $M$ is connected,
$\Scal_+$ is constant almost everywhere, and hence everywhere on $M$ by continuity.

Suppose that $\Scal_+\equiv c$ for some $c>0$.
Then $\Scal\equiv c$ on $M$, and hence $\Delta\Scal=0$.
By \eqref{eq:trace-equation}, we obtain
$\Rm\equiv0$, which contradicts $\Scal\equiv c>0$.
Therefore $\Scal_+\equiv0$, and hence
\begin{equation*}
  \Scal\leq0
\end{equation*}
everywhere on $M$.
Finally, suppose that $\Scal(p)=0$ at some
point $p\in M$. Since $\Scal\leq0$, the function
$\Scal$ attains its global maximum 0 at $p$.
Since $\Delta\Scal\geq0$, the strong maximum principle implies that  $\Scal\equiv0$ on $M$;  see \cite[Theorem 3.5]{GilbargTrudinger1983}, applied in local coordinates, or \cite[Theorem~4.2]{GoffiPediconi} for a formulation on Riemannian manifolds. It then follows from \eqref{eq:trace-equation}
that $\Rm\equiv0$, so $(M,g)$ is flat.

Therefore, if $(M,g)$ is not flat, then
$\Scal<0$ everywhere on $M$.
\end{proof}

We next establish a further flatness result under the assumption that $\Delta\Scal$ is constant, without requiring finite curvature energy.
\begin{proposition}\label{thm:diff_tau}
Let $(M^n,g)$, $n>4$, be a complete Riemannian manifold.
Suppose that $g$ is a critical metric for $\mathcal A$ on  $\mathcal M(M)$. If
$\nabla\Delta\Scal=0$, then $(M,g)$ is flat.

\end{proposition}
\begin{proof}
{{By \eqref{eq:trace-equation} and $n>4$, we have
\begin{equation*}\label{eq:tau-subharmonic}
  \Delta\Scal\geq0.
\end{equation*} }}
Since $\nabla\Delta\Scal=0$ and $M$ is connected,
$\Delta\Scal$ is constant. {{It then follows from \eqref{eq:trace-equation}
that $|\Rm|^2$ is constant.}} 
{{Suppose that $|\Rm|^2>0$. Then, again by
\eqref{eq:trace-equation},
\begin{equation}\label{eq:tau-positive}
  \Delta\Scal > 0.
\end{equation} 
By Lemma~\ref{lem:tau-L2},
$$
  \Scal^2
  \leq n|\Ric|^2
  \leq \frac{n(n-1)}{2}|\Rm|^2,
$$
so $\Scal$ is bounded above. Moreover,
$$
  |\Ric|^2\leq\frac{n-1}{2}|\Rm|^2,
$$
and therefore the Ricci curvature is bounded from below.}} Indeed, for every unit
vector $X$,
$$
  \Ric(X,X)
  \geq-|\Ric|
  \geq-\sqrt{\frac{n-1}{2}}\,|\Rm|.
$$

If $M$ is compact,  there exists a point $p\in M$ at which $\Scal$
has a maximum.
The maximum principle gives
$$
  \Delta\Scal(p)\leq0,
$$
{{which contradicts \eqref{eq:tau-positive}.}}

Suppose that $M$ is noncompact. Since $(M,g)$ is complete, $\Scal$ is bounded above, and the Ricci curvature
is bounded from below, {{the Omori--Yau maximum principle
\cite{Omori1967,Yau1975} applies to $\Scal$.}} Hence there exists a sequence
$\{p_k\}\subset M$ such that
$$
  \Scal(p_k)\to\sup_M\Scal,
  \qquad
  \Delta\Scal(p_k)<\frac1k.
$$
However, $\Delta\Scal$ is a positive constant, say $b>0$.
Thus
$$
  0<b=\Delta\Scal(p_k)<\frac1k
$$
for every $k$, which is impossible for sufficiently large $k$.
Therefore $|\Rm|^2\equiv0$, and hence $(M,g)$ is flat.
\end{proof}

Theorem~\ref{thm:theoremA} reduces
Question~\ref{ques:finite-energy-rigidity}
to the case in which the scalar curvature is
strictly negative everywhere.
In the next section, we derive differential identities
for weakly Einstein metrics that are critical
for $\mathcal A$ and obtain further rigidity results
under additional assumptions.

\section{Differential Identities for Weakly Einstein Critical Metrics}\label{sec:6_identities}

In this section, we derive differential identities for weakly
Einstein metrics that are critical for $\mathcal A$ on
$\mathcal M(M)$. We include dimension four to compare these identities with the known four-dimensional results. Let $(M^n,g)$, $n\geq4$, be a
Riemannian manifold satisfying
\begin{equation*}\label{eq:we_Acritical}
  \check R=\frac{|\Rm|^2}{n}g,
  \qquad
  \grad\A=0.
\end{equation*}
We define the $(0,3)$-tensor $D$ by \begin{equation*}\label{eq:Ddef}
  D_{ijk}:=\nabla_k\Ric_{ij}-\nabla_j\Ric_{ik}.
\end{equation*}
The differential identities derived below relate the Cotton tensor, the scalar curvature, and the tensor $D$. These identities hold without assuming completeness or finite curvature energy.

By definition, $D_{ijk}=-D_{ikj}$.
With the curvature convention (2.1), we have
\begin{equation}\label{eq:D-Bianchi}
  D_{ijk}=\nabla^aR_{aijk}.
\end{equation}
By the definition of the Cotton tensor,
\begin{equation}\label{eq:D-C}
  D_{ijk}
  =
  C_{ijk}
  +
  \frac{1}{2(n-1)}
  \left(
    \nabla_k\Scal\,g_{ij}
    -
    \nabla_j\Scal\,g_{ik}
  \right).
\end{equation}

We first differentiate the weakly Einstein equation.

\begin{proposition}\label{prop:RD-identity}
Let $(M^n,g)$, $n\ge4$, be a weakly Einstein
Riemannian manifold. Then
\begin{equation}\label{eq:R-D}
  R_j{}^{abc}D_{abc}
  =-\frac{n-4}{4n}\nabla_j|\Rm|^2.
\end{equation}
If, in addition, $g$ is critical for $\mathcal A$
on $\mathcal M(M)$, then
\begin{equation}\label{eq:grad-Delta-D}
  \nabla_j\Delta\Scal+nR_j{}^{abc}D_{abc}=0.
\end{equation}
\end{proposition}

\begin{proof}
Taking the covariant divergence of the weakly Einstein equation
$$
  R_{iabc}R_j{}^{abc}=\frac{|\Rm|^2}{n}g_{ij}
$$
and using \eqref{eq:D-Bianchi}, we obtain
\begin{equation}\label{eq:diff-WE}
  R_j{}^{abc}D_{abc}
  +R_{iabc}\nabla^iR_j{}^{abc}
  =\frac1n\nabla_j|\Rm|^2.
\end{equation}
The second Bianchi identity and the symmetries of the curvature tensor
give
\begin{align*}
  0
  &={}
  R_{iabc}\nabla^iR_j{}^{abc}
  +R_{iabc}\nabla_jR^{aibc}
  +R_{iabc}\nabla^aR^i{}_j{}^{bc} \\
  &={}
  2R_{iabc}\nabla^iR_j{}^{abc}
  -\frac12\nabla_j|\Rm|^2.
\end{align*}
Hence
\begin{equation}\label{eq:R-nabla-R}
  R_{iabc}\nabla^iR_j{}^{abc}
  =\frac14\nabla_j|\Rm|^2.
\end{equation}
Substituting \eqref{eq:R-nabla-R} into \eqref{eq:diff-WE} proves
\eqref{eq:R-D}.

Differentiating \eqref{eq:trace-equation} gives
$$
  \nabla_j\Delta\Scal
  =
  \frac{n-4}{4}\nabla_j|\Rm|^2.
$$

Multiplying \eqref{eq:R-D} by $n$, we also obtain
$$
  nR_j{}^{abc}D_{abc}
  =
  -\frac{n-4}{4}\nabla_j|\Rm|^2.
$$
Adding these two identities proves \eqref{eq:grad-Delta-D}.
Note that \eqref{eq:grad-Delta-D} also holds when $n=4$.

\end{proof}

The next divergence identities provide the link between the weakly
Einstein condition and concircular scalar curvature.

\begin{proposition}\label{prop:divergence-identities}
Let $(M^n,g)$, $n\ge4$, be a weakly Einstein
Riemannian manifold. Suppose that $g$ is a critical
metric for $\mathcal A$ on $\mathcal M(M)$. Then
\begin{equation}\label{eq:div-C}
  \nabla^aC_{ija}
  =\frac{1}{2(n-1)}
  \left(
    (\Hess\Scal)_{ij}
    -\frac{\Delta\Scal}{n}g_{ij}
  \right),
\end{equation}
and
\begin{equation}\label{eq:div-D}
  \nabla^aD_{ija}
  =\frac{1}{2n}\Delta\Scal g_{ij}
  =\frac{n-4}{8n}|\Rm|^2g_{ij}.
\end{equation}
\end{proposition}

\begin{proof}
Combining the Euler--Lagrange equation for $\mathcal A$ with the weakly
Einstein condition gives
\begin{equation}\label{eq:Delta-Ric}
  \Delta\Ric_{ij}
  =
  \frac12(\Hess\Scal)_{ij}
  +\Ric_{ia}\Ric_j{}^a
  -\Ric^{ab}R_{iabj}
  +\frac{1}{2n}(\Delta\Scal)g_{ij}.
\end{equation}
Indeed, the last term follows from \eqref{eq:trace-equation}.

Taking the divergence of the Cotton tensor gives
\begin{equation}\label{eq:div-C-expand}
\begin{split}
  \nabla^aC_{ija}
  ={}&
  \Delta\Ric_{ij}
  -\nabla^a\nabla_j\Ric_{ia} 
  -\frac{1}{2(n-1)}
  \left(
    \Delta\Scal g_{ij}
    -(\Hess\Scal)_{ij}
  \right).
\end{split}
\end{equation}
By the Ricci identity and the identity
$\nabla^a\Ric_{ia}=\frac12\nabla_i\Scal$,
\begin{equation}\label{eq:commuted-Ricci}
  \nabla^a\nabla_j\Ric_{ia}
  =
  \frac12(\Hess\Scal)_{ij}
  -\Ric^{ab}R_{iabj}
  +\Ric_{ia}\Ric_j{}^a.
\end{equation}
Substituting \eqref{eq:Delta-Ric} and \eqref{eq:commuted-Ricci} into
\eqref{eq:div-C-expand} gives
\eqref{eq:div-C}.

Finally, taking the divergence of \eqref{eq:D-C} and using
\eqref{eq:div-C}, we obtain
\begin{align*}
  \nabla^aD_{ija}
  &={}
  \nabla^aC_{ija}
  +\frac{1}{2(n-1)}
  \left(
     \Delta\Scal g_{ij}
    -(\Hess\Scal)_{ij}
  \right) \\
  &={}
  \frac{1}{2n}(\Delta\Scal)g_{ij}.
\end{align*}
Equation~\eqref{eq:trace-equation} then yields \eqref{eq:div-D}.
\end{proof}

The preceding formulas also give two useful identities involving the
Cotton tensor directly. From \eqref{eq:D-C} and the trace-free
property of the Cotton tensor, we obtain
\begin{equation*}\label{eq:D-norm}
  |D|^2
  =
  |C|^2+\frac{1}{2(n-1)}|\nabla\Scal|^2.
\end{equation*}
Moreover, substituting \eqref{eq:D-C} into
\eqref{eq:grad-Delta-D} gives
\begin{equation*}\label{eq:grad-Delta-C}
  \nabla_j\Delta\Scal
  +nR_j{}^{abc}C_{abc}
  +\frac{n}{n-1}\Ric_{ja}\nabla^a\Scal
  =0.
\end{equation*}

Recall that a smooth function $f$ is  concircular if
$
  \Hess f=\frac{\Delta f}{n}g.
$
We allow constant functions in this definition. Equation~\eqref{eq:div-C}
therefore gives the following immediate consequence.

\begin{corollary}\label{cor:Cotton-concircular}
Let $(M^n,g)$, $n\geq4$, be a weakly Einstein
Riemannian manifold. Suppose that $g$ is a
critical metric for $\mathcal A$ on $\mathcal M(M)$. Then the Cotton tensor is divergence-free if and only if the
scalar curvature is concircular. 
Equivalently, 
\begin{equation*}\label{eq:Cotton-concircular}
  \nabla^aC_{ija}=0
  \quad\text{if and only if}\quad
  \Hess\Scal
  =\frac{\Delta\Scal}{n}g.
\end{equation*}
\end{corollary}

The identity \eqref{eq:div-D} also explains the special role of
dimension four. When $n=4$, it reduces to
$$
  \nabla^aD_{ija}=0.
$$
If $M$ is closed, then the skew-symmetry of $D$ in its last two indices
gives
$$
  |D|^2=2D^{abc}\nabla_c\Ric_{ab}.
$$
Integration by parts therefore yields
\begin{equation*}
\begin{split}
  \int_M|D|^2\,dv_g
  &=2\int_MD^{abc}\nabla_c\Ric_{ab}\,dv_g \\
  &=-2\int_M\Ric^{ab}\nabla^cD_{abc}\,dv_g
  =0.
\end{split}
\end{equation*}
Hence $D\equiv0$, and the Ricci tensor is a Codazzi 
tensor. 
This recovers the Codazzi property of the Ricci
tensor for closed four-dimensional weakly Einstein
metrics that are critical for $\mathcal A$;
see \cite[Theorem~2.10]{EuhParkSekigawa2011}.

For $n>4$, the same divergence formula gives a different conclusion.

\begin{corollary}\label{cor:divD-flat}
Let $(M^n,g)$, $n>4$, be a weakly Einstein
Riemannian manifold. Suppose that $g$ is a
critical metric for $\mathcal A$ on $\mathcal M(M)$.
If
\begin{equation*}
  \nabla^aD_{ija}=0,
\end{equation*}
then $(M,g)$ is flat. In particular, if the Ricci
tensor is a Codazzi tensor, then $(M,g)$ is flat.
\end{corollary}

\begin{proof}
By \eqref{eq:div-D},
$$
  \nabla^aD_{ija}
  =\frac{n-4}{8n}|\Rm|^2g_{ij}.
$$
If $\nabla^aD_{ija}=0$, then $|\Rm|^2=0$ since $n>4$. Hence $\Rm\equiv0$,
and $(M,g)$ is flat. If the Ricci tensor is a Codazzi tensor, then
$D\equiv0$, so the conclusion follows immediately.
\end{proof}

Corollary~\ref{cor:Cotton-concircular} will be
combined with Theorem~\ref{thm:concircular},
proved in the next section, to establish
Corollary~\ref{cor:divC-rigidity-intro}.

\section{Proofs of Theorem~\ref{thm:concircular}
and Corollary~\ref{cor:divC-rigidity-intro}}\label{sec:7_TheoremB}

We begin with an elementary consequence of geodesic completeness that will be used to prove Theorem~\ref{thm:concircular}, which
does not require the weakly Einstein condition.

\begin{proposition}\label{prop:Hessian}
Let $(M,g)$ be complete. Suppose that a smooth function $f$ satisfies
$$
  \Hess f\geq a f^2 g
$$
for some constant $a>0$. Then $f\equiv0$.
\end{proposition}
\begin{proof}
Let $\gamma:\mathbb R\to M$ be a complete unit-speed geodesic and set
$$
  u(t)=f(\gamma(t)).
$$
Then
$$
  u''(t)
  =
  \Hess f(\dot\gamma(t),\dot\gamma(t))
  \geq a u(t)^2
$$
for all $t\in\mathbb R$. In particular, $u''\geq0$, so $u'$ is
nondecreasing.

{{We claim that $u'\equiv0$. Suppose first that
$u'(t_0)>0$ for some $t_0\in\mathbb R$. Since $u'$ is nondecreasing,
$$
  u'(t)\geq u'(t_0)>0
  \qquad\text{for }t\geq t_0
$$
and hence $u(t)\to+\infty$ as $t\to+\infty$. We may therefore choose
$t_1\geq t_0$ such that
$$
  u(t_1)>0.
$$
For $t\geq t_1$, we have $u'(t)>0$. Multiplying
$u''\geq au^2$ by $u'$ and integrating from $t_1$ to $t$ gives
\begin{equation}\label{u-eq}
  \frac12\bigl(u'(t)^2-u'(t_1)^2\bigr)
  \geq
  \frac a3\bigl(u(t)^3-u(t_1)^3\bigr).
  \end{equation}
Since $u(t)\to+\infty$, we may choose $t_2 \ge t_1$ such that $u(t_2)^3 \ge 2u(t_1)^3$. Since $u$ is increasing on $[t_1,\infty)$, we have $u(t)^3\ge 2u(t_1)^3$ for all $t\ge t_2$. The inequality \eqref{u-eq} therefore gives
$$
  u'(t)^2\geq\frac a3u(t)^3
  \qquad\text{for }t\geq t_2.
$$
Thus
$$
  u'(t)\geq\sqrt{\frac a3}\,u(t)^{3/2}
  \qquad\text{for }t\geq t_2.
$$
It follows that
$$
  \frac{d}{dt}u(t)^{-1/2}
  =
  -\frac12u(t)^{-3/2}u'(t)
  \leq
  -\frac12\sqrt{\frac a3}.
$$
Integrating once more yields
$$
  u(t)^{-1/2}
  \leq
  u(t_2)^{-1/2}
  -
  \frac12\sqrt{\frac a3}\,(t-t_2),
$$
which is impossible for all sufficiently large $t$. Hence
$u'(t_0)>0$ cannot occur.

If $u'(t_0)<0$ for some $t_0$, apply the preceding argument to
$v(t)=u(-t)$. Since
$$
  v''(t)\geq av(t)^2
$$
and $v'(-t_0)=-u'(t_0)>0$, we obtain the same contradiction.
Therefore $u'\equiv0$. Hence $u$ is constant, and
$$
  0=u''\geq au^2
$$
implies $u\equiv0$.}}

Since $\gamma$ is arbitrary, it follows that $f\equiv0$ on $M$.
\end{proof}

\begin{proof}[Proof of Theorem~\ref{thm:concircular}]
Since the scalar curvature $\Scal$ is concircular, \eqref{eq:trace-equation} gives
\begin{equation}\label{eq:Hessian-tau}
 \Hess\Scal=\frac{\Delta\Scal}{n}g=\frac{n-4}{4n}|\Rm|^2g.
\end{equation}
By Lemma \ref{lem:tau-L2},
$$
  |\Rm|^2
  \geq
  \frac{2}{n(n-1)}\Scal^2.
$$
Substituting this estimate into \eqref{eq:Hessian-tau}, we obtain
\begin{equation*}\label{eq:Hessianlower}
  \Hess\Scal
  \geq
  \frac{n-4}{2n^2(n-1)}\Scal^2g.
\end{equation*}
Since $n>4$, the coefficient on the right-hand side is positive.
Proposition~\ref{prop:Hessian}, applied to $f=\Scal$, therefore gives
$$
  \Scal\equiv0.
$$
Finally, \eqref{eq:Hessian-tau}
implies $|\Rm|^2\equiv0$. Thus $(M,g)$ is flat.
\end{proof}
We now prove Corollary~\ref{cor:divC-rigidity-intro}.

\begin{proof}[Proof of Corollary~\ref{cor:divC-rigidity-intro}]
Since $g$ is weakly Einstein and critical for $\mathcal{A},$ Corollary~\ref{cor:Cotton-concircular} applies: 
the assumption that the Cotton tensor is divergence-free implies
\begin{equation*}
  \Hess\Scal=\frac{\Delta\Scal}{n}g.
\end{equation*}
Thus the scalar curvature is concircular, and
Theorem~\ref{thm:concircular} implies that
$(M,g)$ is flat.
\end{proof}

In particular, Corollary~\ref{cor:divC-rigidity-intro} yields flatness when the Cotton tensor vanishes
identically.

Together with Corollary~\ref{cor:divD-flat},
Corollary~\ref{cor:divC-rigidity-intro} shows that
additional differential conditions on $D$ or $C$
imply flatness, without requiring noncompactness
or finite curvature energy.
Question~\ref{ques:finite-energy-rigidity}
remains open without such additional conditions.\\

The following complete noncompact examples show that neither
criticality for $\mathcal A$ nor concircular scalar curvature can be
omitted from Theorem~\ref{thm:concircular}. Criticality is understood
in the local variational sense of Section~\ref{sec:3_EL}.

\begin{example}\label{ex:critical-nonconcircular}
There exists a complete noncompact nonflat metric in dimension five
that is critical for $\mathcal A$ on $\mathcal M(M)$ and whose scalar
curvature is not concircular.

\normalfont Let $(\mathbb H^4,h_{-1})$ be the complete hyperbolic space of
sectional curvature $-1$. Consider
\begin{equation}\label{eq:example-warped-metric}
  M=\mathbb R\times\mathbb H^4,\qquad
  g=dt^2+f(t)^2h_{-1},
\end{equation}
where $f$ is determined by
\begin{equation}\label{eq:example-warped-ode}
  f^2f^{(4)}+4ff'f'''+3f(f'')^2
  -7(f')^2f''-3f''=0
\end{equation}
and
\begin{equation}\label{eq:example-warped-initial}
  f(0)=1,\qquad f'(0)=0,\qquad
  f''(0)=\sqrt{\frac32},\qquad f'''(0)=0.
\end{equation}
We verify below that the unique local solution extends to an even
positive smooth function on $\mathbb R$ and that
\eqref{eq:example-warped-metric} has the asserted properties.

Since the coefficient of $f^{(4)}$ in
\eqref{eq:example-warped-ode} is $f^2$, the initial value problem has
a unique smooth local solution. Equation~\eqref{eq:example-warped-ode} is invariant under
$t\mapsto-t$, so uniqueness and
\eqref{eq:example-warped-initial} imply that $f$ is even.

{{Let $\{\bar e_a\}_{a=1}^4$ be a local orthonormal frame on
$(\mathbb H^4,h_{-1})$. Then
\begin{equation*}
  e_0=\partial_t,\qquad e_a=f^{-1}\bar e_a
  \quad (1\le a\le4)
\end{equation*}
form a local orthonormal frame for $g$.
The curvature components, up to the curvature symmetries, are
\begin{equation}\label{eq:example-curvature-components}
\begin{aligned}
  R_{0ab0}
    &=-\frac{f''}{f}\delta_{ab},\\
  R_{abcd}
    &=-\frac{1+(f')^2}{f^2}
      \bigl(\delta_{ad}\delta_{bc}-\delta_{ac}\delta_{bd}\bigr),\\
  R_{0abc}&=0
\end{aligned}
\end{equation} for
$a,b,c,d\in\{1,2,3,4\}$.

Write $G=\grad\A$.
Using \eqref{eq:example-curvature-components}, we compute from
\eqref{eq:gradA} that
\begin{equation*}
  G=G_0\,dt^2+G_1f^2h_{-1},
\end{equation*}
where}}
\begin{equation}\label{eq:example-gradient-radial}
\begin{aligned}
  G_0={}&\frac4{f^4}\Bigl[
    -2f^2(f'')^2+4f^2f'f'''
    +8f(f')^2f'' \\
    &\hspace{38mm}
    +3\bigl(1+(f')^2\bigr)\bigl(1-3(f')^2\bigr)
    \Bigr],
\end{aligned}
\end{equation}
\begin{equation*}\label{eq:example-gradient-horizontal}
\begin{aligned}
  G_1={}&\frac4{f^3}\Bigl[
    f^2f^{(4)}+4ff'f'''+3f(f'')^2
    -7(f')^2f''-3f''
    \Bigr]. 
\end{aligned}
\end{equation*}
Equation~\eqref{eq:example-warped-ode} gives $G_1=0$.
The diffeomorphism invariance of $\mathcal A$ implies
$\nabla^iG_{ij}=0$, which here reads
\begin{equation*}
  G_0'+4\frac{f'}f(G_0-G_1)=0.
\end{equation*}
By \eqref{eq:example-warped-initial} and
\eqref{eq:example-gradient-radial},
\begin{equation*}
  G_0(0)=12-8f''(0)^2=0.
\end{equation*}
Consequently $(f^4G_0)'=0$ and $G_0=0$. Thus $G=0$ wherever the
solution $f$ is defined.

We next prove global existence and completeness. For small $t>0$,
we have $f>1$ and $f'>0$. Set
\begin{equation*}
  x=\log f,\qquad p=(f')^2,\qquad v=\frac{dp}{dx}=2 f f''.
\end{equation*}
The radial equation $G_0=0$ becomes
\begin{equation}\label{eq:example-phase-system}
  p_x=v,\qquad
  2p\,v_x=\frac12v^2-2pv-3(1-2p-3p^2).
\end{equation}
Let $\sigma=\sqrt6$. Expansion at $t=0$, using
\eqref{eq:example-warped-ode}--\eqref{eq:example-warped-initial},
gives
\begin{equation*}
  p\longrightarrow0,\qquad
  v=\sigma+(\sigma-2)p+O(p^2)
  \quad\text{as }t\downarrow0.
\end{equation*}
In particular, $\sigma<v<2p+\sigma$ initially. On the two boundary curves,
\eqref{eq:example-phase-system} gives
\begin{equation*}
  \left.v_x\right|_{v=\sigma}
    =3-\sigma+\frac92p>0,\qquad
  \left.(v-2p)_x\right|_{v=2p+\sigma}
    =3-2\sigma-\frac p2<0.
\end{equation*}
A first-exit argument therefore yields
\begin{equation*}\label{eq:example-phase-bounds}
  \sigma\le v\le2p+\sigma,\qquad
  \sigma x\le p(x)\le\frac {\sigma}{2}(e^{2x}-1).
\end{equation*}
These bounds keep $p$ positive and prevent finite-time blow-up of
$(p,v)$ for $x>0$. Since
\eqref{eq:example-phase-system} is smooth for $p>0$, its solution
extends to every $x>0$. Moreover, the reconstructed function $f$
satisfies $G_0=0$ and $f'>0$. The divergence identity then gives
$G_1=0$, so this reconstruction continues the original local solution.
We also have
\begin{equation*}
  t=\int_0^x\frac{e^u}{\sqrt{p(u)}}\,du
  \ge\sqrt{\frac{2}{\sigma}}\,x\longrightarrow\infty.
\end{equation*}
Thus $f$ is defined for all $t\ge0$; evenness gives its extension
to $\mathbb R$. Since $f\ge1$,
\begin{equation*}
  g\ge dt^2+h_{-1}.
\end{equation*}
The product metric on the right is complete, and hence $g$ is
complete.

Finally, the scalar curvature and its Hessian are
\begin{equation*}
  \Scal=-8\frac{f''}f
       -12\frac{1+(f')^2}{f^2},\qquad
  \Hess\Scal
  =\Scal''\,dt^2+\Scal'\frac{f'}f(g-dt^2).
\end{equation*}
The initial values and \eqref{eq:example-warped-ode} imply
\begin{equation*}
  \Scal'(0)=0,\qquad \Scal''(0)=12.
\end{equation*}
Therefore $\Hess\Scal=12\,dt^2$ at $t=0$, which is not a scalar
multiple of $g$,  so the scalar curvature is not concircular. The metric is nonflat because $R_{0aa0}(0)=-f''(0)\neq 0$. Hence the
concircularity assumption in Theorem~\ref{thm:concircular} cannot be omitted.
\end{example}

\begin{remark}
The metric of Example~\ref{ex:critical-nonconcircular} also shows that the weakly
Einstein assumption in Corollary~\ref{cor:divC-rigidity-intro} cannot be omitted. Indeed, since the
base is one-dimensional and $(\mathbb{H}^{4},h_{-1})$ has constant
sectional curvature, $(M,g)$ is locally conformally flat, so its Cotton
tensor vanishes identically; yet $g$ is complete, nonflat, and critical
for $\mathcal{A}$ on $\mathcal{M}(M)$. Consistently with Corollary \ref{cor:divC-rigidity-intro}, $g$ is not weakly Einstein: at $t=0$,
\eqref{eq:example-warped-initial} and \eqref{eq:example-curvature-components} give $\check R_{00}=8f''(0)^2=12$ and $\check R_{aa}=2f''(0)^2+6=9$
for $1\le a\le 4$.
\end{remark}

\begin{example}\label{ex:noncritical-concircular}
For every $n>4$, there exists a complete noncompact nonflat metric in
dimension $n$ whose scalar curvature is concircular but which is not
critical for $\mathcal{A}$ on $\mathcal{M}(M)$.
For every $n>4$, the product
\begin{equation*}\label{eq:example-noncritical-product}
  (M^n,g)=
  \bigl(\mathbb R\times S^{n-1},\,dt^2+h_1\bigr),
\end{equation*}
where $h_1$ is the unit round metric, is complete and noncompact.
Its scalar curvature is concircular, but $g$ is not critical for
$\mathcal A$ on $\mathcal M(M)$.

\normalfont Indeed, the product curvature gives
\begin{equation*}
  \Ric=(n-2)h_1,\qquad
  \Scal=(n-1)(n-2),\qquad
  |\Rm|^2=2(n-1)(n-2).
\end{equation*}
In particular, $\Scal$ is constant, so
$\Hess\Scal=(\Delta\Scal/n)g=0$. Since $\Ric$ is parallel,
substitution in \eqref{eq:gradA} yields
\begin{equation*}
  \grad\A
  =(n-1)(n-2)\,dt^2+(n-5)(n-2)\,h_1.
\end{equation*}
Its $tt$-component is nonzero. Equivalently, the trace equation
\eqref{eq:trace-equation} would require
\begin{equation*}
  0=\Delta\Scal
    =\frac{(n-4)(n-1)(n-2)}2>0,
\end{equation*}
which is impossible. Hence $g$ is not critical for $\mathcal A$, and the criticality
assumption in Theorem~\ref{thm:concircular} cannot be omitted.
\end{example}

\end{document}